\documentstyle{article}
\newcommand{\C}{{{\rm C}\!\!\!{\rm I}\,}}
\newcommand{\R}{{{\rm I}\!{\rm R}}}

\newcommand{\re}{{\rm Re}}
\newcommand{\im}{{\rm Im}}
\renewcommand{\epsilon}{\varepsilon}
\renewcommand{\phi}{\varphi}
\begin{document}
\begin{large}

\title{Cartwright-type and Bernstein-type theorems \\
for functions analytic in a cone}
\author{Vladimir Logvinenko
\and Alexander Russakovskii}
\date{}
{\maketitle}

\begin{abstract}
Cartwright-type and Bernstein-type theorems, previously known
only for functions of exponential type in $\C^n,$ are extended
to the case of functions of arbitrary order in a cone.
\end{abstract}

\section {Introduction}

We use standard notations of multidimensional complex analysis.

Let $C$ be an open cone in $\C^n$ with vertex at the origin.
By $[\rho, \sigma, C]$ we denote the class of all
functions $f(z),$ holomorphic in $C$ and satisfying the
estimate
$$\limsup_{|z|\to\infty,\ z\in C} \frac{\log
|f(z)|}{|z|^\rho}\leq \sigma,\ |z|^2=|z_1|^2+\ldots+|z_n|^2.$$

For entire functions we write simply $[\rho, \sigma].$ Thus
$[1,\sigma]$ is the class of entire functions of exponential
type not exceeding $\sigma$ in $\C^n.$

By $[\rho, \infty, C)$ we denote $\displaystyle
\bigcup_{\sigma>0} [\rho, \sigma, C].$

In this paper we make systematical use of the possibility of a
"good" approximation of a function analytic in a cone by entire
functions with control of growth. In the case of dimension $1$
such an approximation was constructed by M.V.Keldysh
\cite{Keldysh}. For the case of several variables, the
result showing the possibility of such an approximation is due
to the second author \cite {Russakovskii}. To formulate this
result we introduce some notations.

Let $\omega$ and $\phi$ be plurisubharmonic functions in
$\C^n,$ both possessing the "non - oscillating" property
$$(u)^{[1]}(z) \leq - A (-u)^{[1]}(z) + B,$$
where  by $u^{[r]}(z)$ we denote $\sup \{u(w):\ |z-w|<r\}.$
Assume also that $\phi(z) \geq 0, \  log(1+|z|) = o(\phi(z)),\
|z| \rightarrow \infty .$ For $\epsilon \geq 0$ we denote by
$\Omega_\epsilon$ the set
$$\Omega_\epsilon = \{ z \in \C^n :\ \omega(z) < - \epsilon
\phi(z) \}$$
and suppose that
$$\forall \epsilon_1>\epsilon_2 :\
inf \{|z_1-z_2|:\ z_1 \in \Omega_{\epsilon_1}, z_2 \in \C^n
\setminus \Omega_{\epsilon_2}\}> 0,$$
which is a kind of smoothness condition on $\omega$ and
$\phi.$

\bigskip

{\bf Theorem A}{\cite {Russakovskii}}.
\begin{it}
Let $f(z)$ be an analytic function in $\Omega_0$ satisfying
the estimate $$|f(z)| \leq C_f e^{C_f \phi(z)}, \  z \in
\Omega_0.$$

Then for each $\epsilon > 0$ and each $N \geq 1$ there exists
such an entire function $g(z)$ that $$|f(z) - g(z)| \leq C
e^{-N\phi(z)}, \  z \in \Omega_\epsilon,$$

$$|g(z)| \leq C  e^{C max (N, C_f) \cdot (\frac{2}{\epsilon}
\cdot \omega^+ + \phi ) (z)}, \  z \in \C^n,$$
where $C$ does not depend on $N.$
\end{it}

\bigskip

We first apply Theorem A to prove results concerning
"boundedness" of functions analytic in a cone. It is natural to
call theorems of this kind Cartwright-type theorems.
Results of such kind were known before only for entire
functions of exponential type.

{\bf Definition.} Let $E$ and $F$ be subsets of $\R^n,$ $E$
being measurable. The set $E$ is called
{\it relatively dense with respect to} $F,$ if for some
positive constants $L$ and $\delta$ and every $x\in F$ $$ |E
\bigcap B(x, L)| \geq \delta.$$

\medskip

Here $|A|$ denotes the Lebesque measure of a (measurable) set
$A,$ and $B(x, L)$ is the ball $\{y\in\R^n:\ |x-y|<L\}.$ The
values of $L$ and $\delta$ are called the density parameters.

\medskip

{\bf Definition.} Let $E$ and $F$ be subsets of $\R^n,$ $E$
being measurable. The set $E$ is called
{\it relatively dense of order $\rho$ with respect to} $F,$ if
its image under the map $x_j\mapsto x_j^\rho, \ j=1,\ldots, n,$
is relatively dense with respect to the image of $F.$

\medskip

{\bf Definition.} A set $E\subset \R^n$ is called an $\epsilon
- net$ for a set $F\subset \R^n$ if for every $x\in F$ there
exists such a point $y\in E$ that
$$|x-y| < \epsilon.$$

{\bf Definition.} A set $E\subset \R^n$ is called an $\epsilon
- net$ {\it of order} $\rho$ for a set $F\subset \R^n$ if its
image under the map $x_j\mapsto x_j^\rho, \ j=1,\ldots, n,$ is
an $\epsilon$ - net for the image of $F$ under this map.

Note that $\epsilon$ - nets may be discrete sets.

Given $\eta\in(0,1),$ denote by $C(\eta)$ the cone in the
positive hyperoctant $\R_+^n$ defined by the relation
$$C(\eta) = \left\{ x\in \R_+^n:\ \min_{j=1, \ldots, n}
x_j \geq \eta \max_{j=1,\ldots,n} x_j\right\}.$$

The results on entire functions of exponential type are
formulated as follows:

\bigskip

{\bf Theorem B} \cite {Logvinenko-cone}.
\begin{it}
Let a set $E$ be relatively dense with respect to $\R_+^n.$
Then each entire function $f\in [1,\infty)$ bounded on $E$ is
bounded on $C(\eta)$ for every fixed $\eta\in(0,1).$

Moreover, for each $\sigma\in(0,\infty)$ there exists such a
finite value $\Delta = \Delta (E, \sigma, \eta),$ that for each
entire function $f\in [1, \sigma]$
$$ \sup_{x\in C(\eta)\setminus B(0, R)}
|f(x)| \leq \Delta \cdot \sup_{x\in E} |f(x)|$$
for some $R=R(f, \eta)<\infty.$

\end{it}

\bigskip

{\bf Theorem B$'$} \cite {Logvinenko-cone}.
\begin{it} Let $E$ be an $\epsilon$ - net for $\R_+^n,$ and let
a number $\eta \in [0, 1)$ be given.

Then there exists such a number
$\sigma_0 = \sigma_0 (n, E, \epsilon, \eta)>0$ that for
every $\sigma \in (0, \sigma_0),$ each function
$f\in[1,\sigma]$ bounded on $E$ is bounded on $C(\eta).$

Moreover, there is such a finite
value $\Delta = \Delta (E, \sigma, \epsilon, \eta),$ that
$$\sup_{x\in C(\eta)\setminus B(0, R)}
|f(x)| \leq \Delta \cdot \sup_{x\in E} |f(x)|$$
for some $R=R(f,\eta)<\infty.$

\end{it}

\bigskip

Below we formulate analogues of the above theorems for
functions holomorphic in cones. Note that we will be
able to consider also functions of order $\rho$
different from 1.

For $\rho\geq 1$ put $W_\rho=\{z=(z_1,\ldots, z_n)\in \C^n: \
-\frac{\pi}{2\rho} < \arg z_j < \frac{\pi}{2\rho}, \
j=1,\ldots, n\}.$ Note that $W_1=\C_+^n$  ($\C_+$ stands for
the right halfplane) and that
$W_\rho\bigcap\R^n = \R^n_+$ for each $\rho.$

Our first result is

{\bf Theorem 1.}
\begin{it}
Let $E$ be a relatively dense set of order $\rho\geq 1$ with
respect to $\R_+^n.$

Then for every $\eta\in(0,1)$ each function $f\in[\rho, \infty,
W_\rho),$ which is analytic in a neighborhood of the origin
and bounded on $E$ is bounded on $C(\eta).$

Moreover, for each $\sigma \in (0, \infty)$ there exists such
a finite value $\Delta = \Delta (E, \sigma, \eta),$ that
$$\sup_{x\in C(\eta)\setminus B(0, R)} |f(x)| \leq \Delta
\cdot \sup_{x\in E} |f(x)|$$ for some $R=R(f, \eta)<\infty.$
\end{it}

\bigskip

The corresponding result for $\epsilon$ - nets is

{\bf Theorem 2.}
\begin{it}
Let $E$ be an $\epsilon$ - net of order $\rho\geq 1$ for
$\R^n_+$ and let a number $\eta \in (0, 1)$ be given.

Then there exists such a number
$\sigma_0 = \sigma_0 (n, E, \epsilon, \eta)>0$ that for
every $\sigma \in [0, \sigma_0),$ each function
$f\in[\rho,\sigma, W_\rho]$ which is analytic near the
origin and bounded on $E$ is bounded on $C(\eta).$

Moreover, there is such a finite
value $\Delta = \Delta (E, \sigma, \epsilon, \eta),$ that
$$\sup_{x\in C(\eta)\setminus B(0, R)}
|f(x)| \leq \Delta \cdot \sup_{x\in E} |f(x)|$$
for some $R=R(f,\eta)<\infty.$

\end{it}

\bigskip

{\bf Remark.} Theorem B$'$ and respectively Theorem 2 may be
formally slightly strengthened by assuming the set $E$ to be
an $\epsilon$ - net not for the whole $\R^n_+$ but for its
relatively dense subset.

{\bf Remark.} To the best of our knowledge, Theorems 1 and 2
are new in the case $\rho>1$ even for entire functions.

{\bf Remark.} There are examples showing sharpness (in a
certain sense) of theorems B and B$'$; for instance, it was
shown \cite{Logvinenko-cone} that the results fail to hold if
we do not truncate the cone $C(\eta)$ by the ball $B(0, R),$
and that the value of $R$ cannot be chosen independent of $f\in
[1,\sigma],$ etc. The same examples with obvious
modifications play a similar role for theorems 1 and 2.

\bigskip

Next we mention V.Bernstein-type theorems for entire functions
of finite order. By this we mean results giving conditions on
sets sufficient for calculation of the (radial) indicator.
Remind that the radial indicator of a function $f(z)\in
[\rho,\infty)$ is defined as follows:
$$h_f (z) =\limsup_{w\to z} \limsup_{t\to\infty}\frac{\log
|f(tw)|}{t^\rho}.$$
For the case of dimension $1$ the first $\limsup$
(regularization) may be omitted. We refer to \cite{Ronkin} for
the properties of the radial indicator.

V.Bernstein \cite{Bernstein} was the first to give a sufficient
condition on a set $E$ on a ray which guarantees that $$h_f(1)
= \limsup_{t\to\infty,\ t\in E} \frac{\log |f(t)|}{t^\rho}.$$
The references to the further results in this direction are
given in \cite{Logvinenko-indicator}. We mention below results
of the first author concerning entire functions in $\C^n.$

{\bf Definition.} Let $\epsilon(R),\ R\in \R_+,$ be a function
monotonically decreasing to zero as $R\to\infty.$ A set $E\subset
\R^n$ is called an $\epsilon (R) - net$ for a set $F\subset
\R^n$ if for each $x\in F$ there exists $y\in E$ such that
$$|x-y| \leq \epsilon (|x|).$$

{\bf Definition.} A set $E\subset \R^n$ is called an
$\epsilon(R) - net$ {\it of order} $\rho$ for a set $F\subset
\R^n$ if its image under the map $x_j\mapsto x_j^\rho, \
j=1,\ldots, n,$ is an $\epsilon(R)$ - net for the image of $F$
under this map.

\bigskip

{\bf Theorem C} \cite {Logvinenko-indicator}.
\begin{it}Let a set $E$ be an $\epsilon(R)$ - net of order
$\rho\in (0,\infty)$ for some cone $C(\eta_0).$

Then the relation
$$h_f(\frac{1}{\sqrt n},...,\frac{1}{\sqrt n}) = \lim_{\eta\to
0} \limsup_{|x|\to\infty,\ x\in E\bigcap C(\eta)}
\frac{\log |f(x)|}{|x|^\rho}$$
holds for every function $f\in [\rho, \infty).$

\end{it}

\bigskip

Theorem C yields the following uniqueness result.

\bigskip

{\bf Theorem D} \cite {Logvinenko-indicator}.
\begin{it}Let $E$ be as in theorem C and let
$$\limsup_{|x|\to\infty, \ x\in E} \frac{\log
|f(x)|}{|x|^\rho}=-\infty$$
for some function $f\in [\rho,\infty).$

Then $f(z)\equiv 0.$

\end{it}

\bigskip

The cones $W_\tau$ were defined above for $\tau\geq 1$. Now we
would like to extend the definition to all $\tau>\frac{1}{2}.$
For $\tau\in (\frac{1}{2}, 1)$ define $W_\tau$ to be the same as
$W_{\frac{2\tau}{2\tau -1}}.$

Our theorem 3 below is an analogue of theorem C for functions
holomorphic in cones.

\bigskip

{\bf Theorem 3.}
\begin{it}Let a set $E$ be an $\epsilon(R)$ - net of order
$\rho\geq \frac{1}{2}$ for some cone $C(\eta_0).$

Then the relation
$$h_f(\frac{1}{\sqrt n},\ldots,\frac{1}{\sqrt n}) =
\lim_{\eta\to 0} \limsup_{|x|\to\infty,\ x\in E\bigcap C(\eta)}
\frac{\log |f(x)|}{|x|^\rho}$$
holds for every function $f\in [\rho, \infty, W_\tau), \
\tau>0.$

\end{it}

\bigskip

Note that while the indicator of an entire function of finite
type $\sigma$ is bounded below by $-\sigma$ \cite{Ronkin}, the
(regularized) indicator of a function holomorphic in a cone
needs not to be bounded from below. Hence the corresponding
uniqueness result holds only if the cone $W_\tau,$ in which our
function is defined, is wide enough.

\bigskip

{\bf Theorem 4}
\begin{it}Let $E$ be as in theorem 3 with $\rho>2$ and let
$$\limsup_{|x|\to\infty, \ x\in E} \frac{\log
|f(x)|}{|x|^\rho}=-\infty$$
for some function $f\in [\rho,\infty, W_\tau), \ \tau\leq
\frac{\rho}{2}.$

Then $f(z)\equiv 0.$

\end{it}

\bigskip

\section {Some remarks concerning cones in $\C^n$}

In this paper we will deal mainly with two types of cones in
$\C^n$.

One of them, $W_\tau,$ is defined in the previous section. We
introduce another one.

For $t>0$ denote by $||.||_t$ a norm in $\C^n$ given by
$$||z||_t = \max_{j=1,\ldots, n} \{ |\re z_j|,\ |\im
\frac{z_j}{t}|\}.$$

By $Y_t (\eta), \ \eta\in[0,1),$ we denote the cone in $\C^n$
given by
$$Y_t (\eta)=\left\{ z\in \C_+^n:\ \min_{j=1, \ldots, n} \re
z_j\geq \eta ||z||_t\right\}.$$

Note that for all $t>0$ the intersection of $Y_t (\eta)\bigcap
\R_+^n$ is exactly the real cone $C(\eta).$ Obviously,
$Y_t(0)=\C_+^n.$

The geometry of the cone $Y_t (\eta)$ is very simple. We just
observe that the ray ${\ell}= \{\xi (1,\ldots, 1),
\xi>0\}$ lies on the complex line ${\cal L}=\{z_1=\ldots=z_n\}$
which has the largest intersection with $Y_t (\eta):$
$${\cal L} \bigcap Y_t(\eta)= \{w\in \C: |\arg w|<\arctan
t\eta\}.$$

One easily sees that, given a number $\tau>0,$ it is
possible to choose such $t$ and $\eta$ that
$$Y_{t}(\eta)\subset W_\tau$$ and
$${\cal L} \bigcap Y_t(\eta) = {\cal L} \bigcap W_\tau.$$

We would like to write each of the two types of cones in the
form $\{z\in \C^n:\ u(z)<0\}$ for some plurisubharmonic
function $u(z)$ in $\C^n.$ For $Y_t(\eta)$ we can take $u(z)$
to be of order $1:$ $$u(z)=\max_{j=1,\ldots,n}(-\re z_j) + \eta
||z||_t,$$ while for $W_\tau$ with $\tau\in(1/2, 1)$ one can take
function of any order $\rho\in[\tau, 1):$
$$u(z)=\max_{j=1,\ldots,n}u_j(z),$$
where
$$u_j = \left\{
\begin{array}{ll}
|z_j|^\rho \sin \rho (\arg z_j - \frac{\pi}{2\tau}),&
\mbox{if $|\arg z_j|<\frac{\pi}{2\tau}+\frac{\pi}{2\rho}$;}\\
|z_j|^\rho, & \mbox{if $|\arg z_j|\in
[\frac{\pi}{2\tau}+\frac{\pi}{2\rho},\pi].$}\\
\end{array}
\right.
$$

\bigskip

\section {Proof of Theorems 1 and 2}

{\bf Proof of Theorem 1.} The idea of the proof is to
approximate the function $f(z)$ by an entire function $g(z)$
with the help of theorem A, apply theorem B to $g(z)$ and
derive the required estimates for $f(z).$

First we note that theorem B may be reformulated for an
arbitrary cone $C\subset \R^n,$ since it is always possible to
find such an automorphism $\psi:\C^n\to\C^n,$ that $\psi
(\R_+^n) \subset C,$ and to consider $f(\psi(z))$ instead of
$f(z)$ which results in an obvious recalculation of all
coefficients and does not affect the order of the holomorphic
function.

Due to the possibility of the transformation $z_j\mapsto
z_j^\rho,\ j=1,\ldots, n$ which takes cones (with vertex in
the origin) into cones, particularly, $C(\eta)$ into
$C(\eta^{1/\rho}),$ relatively dense sets of order $\rho$ into
relatively dense sets of order $1$ and holomorphic functions of
order $\rho$ in $W_\rho$ into holomorphic functions of order
$1$ in $W_1=\C_+^n$ it is enough to assume $\rho=1$ in what
follows.

Next we define functions $\omega$ and $\phi$ in the following
way. Put $$\omega = \max_{j=1,\ldots,n} (-\re
z_j), \ \phi = \max (\delta, ||z||_t), $$ where $t>0$ is
arbitrary, and we choose $\displaystyle
\delta=\log\frac{1}{\sup_{x\in E} |f(x)|}.$

Then the set $\Omega_0 = \{ z:\ \omega(z)<0\}$ is exactly
$\C_+^n,$ and for $\epsilon \in (0,1)$ the set
$\Omega_\epsilon= \{z:\ \omega(z) < - \epsilon
\phi(z)\}$ (which is $Y_t (\epsilon)$ without some
neighborhood of the vertex) has the property $\Omega_\epsilon
\bigcap \R_+^n \supset C(\epsilon)\setminus B(0, \epsilon
\delta).$ It is clear that the conditions of theorem A are
satisfied.

Let $f(z)$ be a given function from the class $[1,\sigma, W_1]$.
By theorem A, there exists such an entire function $g(z)$ of
exponential type $\leq K\sigma$ with $K=K(n, \epsilon)$ not
depending on $f\in [1,\sigma]_+$ that $|f(z)-g(z)|\leq
e^{-\phi}, \ z\in \Omega_\epsilon.$ Our choice of $\delta$
implies that
$$ \sup_{x\in E\bigcap (C(\epsilon)\setminus B (0,
\epsilon\delta)} |g(x)| \leq 2 \sup_{x\in E} |f(x)|.$$

The function $f^* (z) = g(z-\delta)$ is an entire
function belonging to the class $[1, K\sigma]$ bounded by
$\displaystyle 2\sup_{x\in E} |f(x)|$ on a set
$E^*=\{z+\delta:\ z\in E\}$ which is relatively dense with
respect to the cone $C(\epsilon).$ By our remark above we can
apply theorem B to $f^*.$ According to this theorem, for each
$\eta\in(0,\epsilon)$ there exist such positive numbers
$\Delta$ and $R=R(f)$ that $$\sup_{x\in C(\eta)\setminus
B(0,R)} |f^*(x)| \leq \Delta\cdot\sup_{x\in E^*} |f^*(x)|\leq
2\Delta\cdot\sup_{x\in E} |f(x)|.$$

For each $\eta_1\in (0, \eta)$ there exists such a number
$R_1>R$ that $$C(\eta_1)\setminus B(0,R_1) \subset
\{z+\delta: z\in C(\eta)\setminus B(0, R)\}.$$

Since $$|f(x)-f^*(x+\delta)|\leq
\sup_{t\in E} |f(t)|$$ for $x\in C(\eta_1)\setminus B(0, R_1),$
and the numbers $\epsilon\in (0,1),\ \eta \in (0,\epsilon),
\eta_1 \in (0, \eta)$ were arbitrarily chosen, we obtain the
required estimate.

The theorem is proved.

\bigskip

{\bf Proof of Theorem 2.} The proof repeats the proof of the
previous theorem with the only difference that theorem B$'$ is
applied instead of theorem B. The corresponding value of
$\sigma_0$ in theorem 2 differs from that in theorem B$'$ by
the factor $\frac{\max (\eta/2, 1)}{C}$ where $C$ is the
constant from theorem A.

\bigskip

\section {Proof of Theorems 3 and 4}

{\bf Proof of Theorem 3.}  Given a function $F(z)$ analytic in
$W(\tau)$ and of order $\rho,$ denote
$h_F(\frac{1}{\sqrt n},\ldots,\frac{1}{\sqrt n})$ by $H_F$
and let $H_F (E)$ be the corresponding limit calculated over
the set $E.$ It is obvious that $$H_F (E) \leq H_F.$$ We need
to prove the converse.

The way to do it is to use theorem A to find an
entire function $g\in [\rho,\infty)$ with the properties
$$H_g=H_f$$ and $$H_g(E)=H_f(E)$$ and use theorem $C$ for entire
functions to prove that $$H_g(E)=H_g,$$ which yields the desired
relation.

Assume that $\rho\geq 1$ first. Then one can
take $\rho=1$ by the same arguments as before.
We consider two cases. Assume first that $H_f(E)>-\infty.$ Since
multiplication of our function by $e^{A(z_1+\ldots+z_n)}$ does
not affect the investigated property of the set $E,$ we can
always assume that $H_f(E)> 0.$ Hence an entire function $g(z)$
uniformly approximating $f(z)$ in a cone containing the ray
${\ell} = \{(t,\ldots,t),\ t>0 \}$ will have the same values
of $H_g(E)$ and $H_g$ as $f.$ Thus it is enough to construct
such a function.

Take $\displaystyle N=1, \
\omega(z) = \max_{j=1,\ldots, n} (-\re\ z_j)+\eta'||z||_t,$
for such $t$ and $\eta'$ that the cone $\Omega_0 =
\{z:\ \omega(z)<0 \}$ is contained in $W_\tau,$ and take
$\phi(z) = \max (\delta, ||z||_t).$ For $\epsilon \in (0,
1-\eta')$ the set $\Omega_\epsilon = \{z:\ \omega(z) < -
\epsilon\phi(z)\}$ (which is $Y_t (\eta'+\epsilon)$ without a
neighborhood of the vertex) has the property $$\Omega_\epsilon
\bigcap \R_+^n \supset C(\eta'+\epsilon)\setminus B(0,
\epsilon\delta).$$ Applying theorem A, we are done.

Now consider the case $H_f(E)=-\infty.$ We need to prove that
$H_f=-\infty.$ We choose $\omega$ and $\phi$ to be the same as
above. Fix some $\epsilon\in (0,1-\eta)$ and denote by $g_N$
the entire function of finite type corresponding to the choice
of $N\geq 1$ in theorem A. Note that for any such function
$$H_{g_N}(E)=H_{g_N}.$$

The entire function $g_N(z)$ satisfies $|f(z) -
g(z)|<e^{-N||z||_t}$ on $\Omega_\epsilon,$ in particular, on
$C(\epsilon')$ for $|x|$ large enough.
We have
\begin{eqnarray*}
H_{g_N}& \leq &
\lim_{\eta\to 0}\limsup_{|x|\to \infty,\ x\in E\bigcap C(\eta)}
\frac{\log (|f(x)|+|g_N(x)-f(x)|)}{|x|}\\
&\leq & \max (-N, H_f(E))\\
&= & -N,\\
\end{eqnarray*}
and
\begin{eqnarray*}
H_{f} &\leq &
\lim_{\eta\to 0}\limsup_{|x|\to \infty, x\in C(\eta)}
\frac{\log (|g_N(x)|+|g_N(x)-f(x)|)}{|x|} \\
&\leq& \max (-N, H_{g_N}(E))\\
& = & \max(-N, H_{g_N})\\
& = & -N.
\end{eqnarray*}

Since $N$ was arbitrary, we conclude that $H_f=-\infty.$ The
theorem is proved in this case.

The proof in the case $\rho\in(\frac{1}{2}, 1)$ follows the
same scheme as above with the only difference that we do not
pass over to the order $1.$ To apply theorem A, for $\omega
(z)$ we take the function $u(z)$ mentioned in the end of
section 2, and set $\displaystyle \phi(z)=\max_{j=1,\ldots,n}
(|z_j|^\rho, \delta).$

\bigskip

{\bf Remark.} In the case $\rho\geq 1$ it is also possible to
give another proof of Theorem 3 based on Theorem 2 (and thus
using Theorem A indirectly).

\bigskip

{\bf Proof of Theorem 4.} Since the conditions of Theorem 4
imply that the rays ${\ell} =
\{(t,\ldots,t),\ t>0\}$ and $e^{i\frac{\pi}{\rho}}{\ell}$ both
belong to the cone where our function is holomorphic, the
result follows from Theorem 3 and the properties of the
indicator \cite{Ronkin} Ch. 3, \S 5.

\bigskip

\section {Acknowledgments}

The authors would like to thank S.Yu.Favorov, A.Yu.Rashkovskii,
L.I.Ronkin and M.L.Sodin for helpful discussions.

\begin{tabular}{ll}
Vladimir Logvinenko \quad & Alexander Russakovskii\\
6644 Rosemead Boulevard & Theory of Functions
Department\\
Apt. 20 & Mathematical Division \\
San Gabriel & Institute for Low Temperature Physics\\
CA 91775 & 47 Lenin Avenue\\
USA & 310164 Kharkov\\
& Ukraine\\
\end{tabular}

\end{large}

\bigskip

\begin{thebibliography}{MMM}

\bibitem[Be]{Bernstein} V.Bernstein. Sur les propri\'et\'es
caract\`eristiques des indicatrices de croissance. C.R., 202,
1936, 108 - 110.

\bibitem[Ke] {Keldysh} M.V.Keldysh. Sur
l'approximation des fonctions holomorphes par les fonctions
enti\`eres. Dokl. AN SSSR, 47, 1946,  239 - 241.

\bibitem[L1]{Logvinenko-cone} V.N.Logvinenko. Boundedness
conditions for entire functions of exponential type inside the
hyperoctant $\R_+^n.$ Math. USSR - Izv., 34, 1990, 663 -
676.

\bibitem[L2]{Logvinenko-indicator} V.N.Logvinenko.
Radial indicator of an entire function may be
calculated over a discrete subset of a subspace of
small dimension. Preprint.

\bibitem [Ro]{Ronkin} L.I.Ronkin. Introduction to the theory of
entire functions of several variables. Transl. of Math.
Monographs, vol. 44, Amer.Math.Soc., Providence, R.I., 1974.

\bibitem[Ru]{Russakovskii} A.M.Russakovskii.
Approximation by entire functions on unbounded domains in
$\C^n.$ J.Approx.  Theory, 74, 1993, n3, 353 - 358.

\end{thebibliography}
\end{document}